\documentclass[12pt]{article}
\usepackage[a4paper, total={7in, 10in}]{geometry}
\usepackage[parfill]{parskip}
\usepackage{physics, tensor, float, subcaption}
\usepackage{graphicx}
\graphicspath{ {Plots/} }
\usepackage{jhep-mod}
\usepackage{bm}
\usepackage{soul}
\usepackage{amssymb,amsmath,amsthm}
\usepackage{mathrsfs}
\usepackage[utf8]{inputenc}
\usepackage{enumerate}
\usepackage{bigints}
\usepackage{xcolor}
\usepackage{appendix}
\usepackage{graphicx}
\usepackage{float}
\usepackage{tikz}
\usepackage{setspace}
\usepackage{cancel}
\usepackage{doi}
\definecolor{purple}{rgb}{1,0,1}
\definecolor{lime}{HTML}{A6CE39} % needs xcolor
\newcommand{\blue}[1]{{\color{blue} #1}}

\definecolor{lime}{HTML}{A6CE39}
\newcommand{\orcidicon}{%
	\begin{tikzpicture}
	\draw[lime, fill=lime] (0,0) 
		circle [radius=0.16] 
		node[white] {{\fontfamily{qag}\selectfont \tiny ID}};
	\draw[white, fill=white] (-0.0625,0.095) 
		circle [radius=0.007];
	\end{tikzpicture}
	\hspace{-5mm}
}
\newcommand\orcidMatt{{\href{https://orcid.org/0000-0003-1088-6485}{\orcidicon}}}
\renewcommand{\O}{\mathcal{O}}
\begin{document}
%========================================================
%========================================================
%========================================================

\title{\null\vspace{-50pt}
\leftline{\blue{Effective short intervals containing primes}}
}

%========================================================
%========================================================
%========================================================
\author{
\Large
Matt Visser\!\orcidMatt\!
}
%========================================================
%========================================================
%========================================================
%========================================================
\affiliation{School of Mathematics and Statistics, Victoria University of Wellington, \\
\null\qquad PO Box 600, Wellington 6140, New Zealand.}
%========================================================
%========================================================
\emailAdd{matt.visser@sms.vuw.ac.nz}
%========================================================
%========================================================
\def\theta{\vartheta}
\def\O{{\mathcal{O}}}
\def\Li{{\mathrm{Li}}}

\abstract{
\vspace{1em}

95 years ago
Hoheisel proved the existence of primes in the sub-linear interval
\[
\left[x, x+x^{1-{1\over 33000}}\right] \qquad \hbox{for $x$ sufficiently large}.
\]
This was improved by Heilbronn in 1933,  proving existence of primes in the sub-linear interval
\[
\left[x, x+x^{1-{1\over 250}}\right] \qquad \hbox{for $x$ sufficiently large}.
\]
More recently Baker, Harman, Pintz proved existence of primes in the sub-linear interval
\[
\left[x, x+ x^{1-{19\over 40}}\right] \qquad \hbox{for $x$ sufficiently large}.
\]
In the present article I will, to the extent possible, make some of these statements effective. 
Specifically, among other things,  I shall show that
\[
\forall n \geq 4, \qquad\forall x \geq \exp(3\exp(33)), \qquad \hbox{there are primes in the interval} \left[x, x+ x^{1-{1\over n}}\right];
\]
\[
\forall n \geq 91, \qquad\forall x \geq [90^{90}]^{n/(n-90)} , \qquad 
\hbox{there are primes in the interval} \left[x, x+ x^{1-{1\over n}}\right].
\]
Furthermore
\[
\forall n \geq 106, \qquad\forall x \geq 1, \qquad 
\hbox{there are primes in the interval} \left[x, x+ x^{1-{1\over n}}\right].
\]
In particular this last observation makes both the Hoheisel  and Heilbronn results fully explicit and effective.
This (relatively) specific observation can be extended and generalized in various ways.

\bigskip

\bigskip
\noindent
{\sc E-print: 2508.18786 [math.NT]}

\bigskip
\noindent
{\sc Date:} 28 August 2025; 12 July 2026; \LaTeX-ed \today

\bigskip
\noindent{\sc Keywords}: Primes in short intervals; effective bounds; Hoheisel; Heilbronn.

\bigskip
\noindent{\sc Changes:} Minor typos fixed; minor improvements in presentation; added appendix.

}

%========================================================
\maketitle
%========================================================
\def\tr{{\mathrm{tr}}}
\def\diag{{\mathrm{diag}}}
\def\cof{{\mathrm{cof}}}
\def\pdet{{\mathrm{pdet}}}
\def\d{{\mathrm{d}}}
\parindent0pt
\parskip7pt
\def\Kerr{{\scriptscriptstyle{\mathrm{Kerr}}}}
\def\eos{{\scriptscriptstyle{\mathrm{eos}}}}

\def\Z{{\mathbb{Z}}}
%================================================
\section{Introduction}
%================================================

\enlargethispage{50pt}
The investigation  of primes in sub-linear intervals of the form $\big[x, x+ x^{1-{1\over n}}\big]$ has, by now, almost a century of history. Unfortunately almost everything we know is ``ineffective'', meaning ``true for sufficiently large numbers'', but with no real understanding of just how large ``sufficiently large'' has to be.
Examples of this behaviour include Hoheisel's 95 year old 
proof~\cite{Hoheisel} of the existence of primes in the sub-linear interval
\begin{equation}
\left[x, x+ x^{1-{1\over 33000}}\right] \qquad \hbox{for $x$ sufficiently large},
\end{equation}
and Heilbronn's 1933 improvement~\cite{Heilbronn},  proving existence of primes in the  smaller sub-linear interval
\begin{equation}
\left[x, x+ x^{1-{1\over 250}}\right] \qquad \hbox{for $x$ sufficiently large}.
\end{equation}
More recently Baker, Harman, and Pintz~\cite{BHP} proved existence of primes in the even smaller sub-linear interval
\begin{equation}
\left[x, x+ x^{1-{19\over 40}}\right] \qquad \hbox{for $x$ sufficiently large}.
\end{equation}
Baker, Harman, and Pintz comment that ``with enough work this result could be made explicit'', but as yet there has been no progress along these lines.

By combining two known effective bounds on primes between powers, with some information on maximal prime gaps, and some recent improvements on primes in effective linear intervals,  I will make some progress on this issue.

%\clearpage
%================================================
\section{Primes between powers}
%================================================

To extract effective bounds on primes in the interval $\big[x, x+x^{1-{1\over n}}\big]$ I shall work backwards from known effective bounds on primes between powers, ostensibly in intervals of the form $\left[i^m, (i+1)^{m}\right]$. However, statements regarding primes between powers 
$\left[i^m, (i+1)^{m}\right]$ are easily converted,  using the identity
\begin{equation}
x = i^m;  \qquad 
(i+1)^m = i^m \left(1+{1\over i}\right)^m = i^m + m i^{m-1}+\dots = 
x + m x^{1-{1\over m}} + \dots,
\end{equation}
into statements regarding intervals of the 
form $\big[ x , x + m \; x^{1-{1\over m}}\big]$, and it is intervals of this form that will be most directly useful below. 

%================================================
\subsection{Primes between 90th powers}
%================================================

Cully-Hugill has shown~\cite{Cully-Hugill}
\begin{equation}
\forall i \in \Z^+ \qquad (i^{155}, [i+1]^{155}) \hbox{\quad contains at least one prime.}
\end{equation}
This has subsequently been updated twice, by Cully-Hugill  and Johnston.\\ 
First, in reference~\cite[Theorem 1.3]{CHJ1}
\begin{equation}
\forall i \in \Z^+ \qquad (i^{140}, [i+1]^{140}) \hbox{\quad contains at least one prime.}
\end{equation}
Second, in reference~\cite[Theorem 1.4]{CHJ2}
\begin{equation}
\forall i \in \Z^+ \qquad (i^{90}, [i+1]^{90}) \hbox{\quad contains at least one prime.}
\end{equation}

Inspection of the proof of~\cite[Theorem 1.4]{CHJ2} shows that what was actually proved is the slightly stronger statement:
\begin{equation}
\forall x\geq 1 \quad \hbox{the interval} \quad [x,  x+ 90 \; x^{1-{1\over90}}] \hbox{\quad contains at least one prime.}
\end{equation}
Now let $n> 90$ and compare the two intervals
\begin{equation}
[x,  x+ x^{1-{1\over n}}],
\qquad
[x,  x+ 90 \; x^{1-{1\over90}}].
\end{equation}

\clearpage
For large $x$  the width of the Hoheisel interval will overtake and subsume that of the Cully--Hugill+Johnston interval, and this will happen when
\begin{equation}
x^{-{1\over n}} = 90 x^{-{1\over90}};  \qquad x^{{1\over 90}-{1\over n}} = 90; 
\qquad x^{{n-90\over 90 n}} = 90;
\qquad x = [90^{90}]^{n\over n-90}.
\end{equation}
Consequently we have the effective result
\begin{equation}
\forall n > 90, \quad \forall x > [90^{90}]^{n\over n-90}, \qquad 
[x,  x+ \; x^{1-{1\over n}}] \hbox{\quad contains at least one prime.}
\end{equation}
This is our first key result, which we shall soon combine with other results to be derived below. To set the scale, note $90^{90}\approx\exp(405)\approx 7.62\times10^{175}$.
Note that by construction this particular technique provides no information for $n\leq 90$.  A few key values are reported in Table 1.

\begin{table}[!htbp]
\caption{A few key instances of the Cully-Hugill+Johnston lower bound.}
\begin{center}
\begin{tabular}{|c|c|c|}
\hline\hline
$n$ & $[90^{90}]^{n\over n-90}$ & Notes \\
\hline\hline
90 & $\infty$ & C-H+J exponent\\
91 & $1.762594084\cdot10^{16005}$ & --- \\
100 & $6.580493968\cdot10^{1758}$ & ---\\
105 & $1.488610050\cdot10^{1231}$ & --- \\
106 & $1.648527716\cdot10^{1165}$ & --- \\
125 & $1.410516056\cdot10^{628} $ &--- \\
250 & $6.536614755\cdot10^{274}$ & Heilbronn exponent\\
500 & $3.090519754\cdot10^{214}$ & ---\\
1000 & $1.891174930\cdot10^{193}$ & ---\\
33000 & $2.305766611\cdot10^{176}$ & Hoheisel exponent\\
$\infty$ &  $7.617734805\cdot10^{175}$ & --- \\
\hline\hline
\end{tabular}
\end{center}
\label{T:1}
\end{table}%

This already makes both the Heilbronn and Hoheisel results effective, but we shall soon be able to say much more.

%================================================
\subsection{Primes between cubes}
%================================================

Dudek has shown~\cite{Dudek} (slightly updated by Cully-Hugill~\cite{Cully-Hugill})\begin{equation}
\forall i \geq \exp(\exp(33))  \qquad (i^{3}, [i+1]^{3}) \hbox{\quad contains at least one prime.}
\end{equation}
Dudek~\cite{Dudek} used the specific constant $\exp(\exp(33.217))$, whereas Cully-Hugill~\cite{Cully-Hugill} slightly improved this to $\exp(\exp(32.892))$. In the interests of clarity of presentation I will round this to $\exp(\exp(33))$.

\clearpage
Inspection of the proof shows that what was actually proved is the slightly stronger statement that:
\begin{equation}
\forall x\geq  \exp(3\exp(33))   \qquad [x,  x+ 3 \; x^{1-{1\over3}}] \hbox{\quad contains at least one prime.}
\end{equation}
Now let $n> 3$ and compare the two intervals
\begin{equation}
[x,  x+ x^{1-{1\over n}}]
\qquad
[x,  x+ 3 \; x^{1-{1\over3}}]
\end{equation}
Eventually, in fact quite rapidly, the width of the Hoheisel interval will overtake and subsume the Dudek interval, and this will happen when
\begin{equation}
x^{-{1\over n}} = 3 x^{-{1\over3}};  \qquad x^{{1\over 3}-{1\over n}} = 3; \qquad x^{{n-3\over 3 n}} = 3;
\qquad x = [27]^{n\over n-3}.
\end{equation}
But we now also need to keep track of the fact that the Dudek interval is only guaranteed to contain a prime for $x> \exp(3\exp(33))$.

Consequently we have the effective result that
\begin{eqnarray}
&&\forall n > 3, \quad \forall x > \max\left\{ \exp(3\exp(33)), \; [27]^{n\over n-3}\right\}, \quad 
\nonumber\\
&&
\qquad
\hbox{the interval } [x,  x+ \; x^{1-{1\over n}}] \hbox{ contains at least one prime.}
\end{eqnarray}
This is our second key result, one  that we shall combine with the result derived above, before ultimately combining it with other results to be derived below.
Note that by construction this particular technique provides no information for $n\leq 3$. 

We can also use this Dudek discussion to slightly improve the Cully-Hugill+Johnston  bound presented above:
\begin{eqnarray}
&&\forall n > 90, \quad \forall x > \min\left\{ \exp(3\exp(33)), \; [90^{90}]^{n\over n-90}\right\}, \quad 
\nonumber\\
&&
\qquad
\hbox{the interval } [x,  x+ \; x^{1-{1\over n}}] \hbox{ contains at least one prime.}
\end{eqnarray}

%================================================
\subsection{Combining these two lower bounds}
%================================================

\enlargethispage{40pt}
%Note that
%\begin{equation}
%\exp(3\exp(33)) =[27]^{n\over n-3}
%\end{equation}
%has two real roots.
The effective bound coming from the Dudek result has an implied crossover point $n_1$, located at 
\begin{equation}
 \exp(3\exp(33)) =[27]^{n_1\over n_1-3}; \qquad  3\exp(33) = {n_1\over n_1-3} \ln 27; 
 \end{equation}
 implying
 \begin{equation}
 \qquad n_1 = {3\over 1 - \exp(-33)\ln 3} \approx 3 + 1.5355 \times 10^{-14}.
\end{equation}
There is a second  crossover point $n_2$,  coming from combining the Cully-Hugill+Johnston result with Dudek's result,  which is located at
\begin{equation}
 \exp(3\exp(33)) =[90^{90}]^{n_2\over n_2-90}; \qquad  3\exp(33) = {n_2\over n_2-90} 90\ln 90; 
 \end{equation}
 implying
 \begin{equation}
 \qquad n_2 = {90\over 1 - \exp(-33) 30\ln 90} \approx  90+  5.6603 \times 10^{-11}.
\end{equation}

%\clearpage
Then we can summarize the results of the previous two subsections as follows:
\begin{equation}
\forall n \in(3,n_1], \quad \forall x > [27]^{n\over n-3}, \qquad 
[x,  x+ \; x^{1-{1\over n}}] \hbox{\quad contains at least one prime.}
\end{equation}
\begin{equation}
\forall n \in[n_1,n_2], \quad \forall x > \exp(3\exp(33)), \qquad 
[x,  x+ \; x^{1-{1\over n}}] \hbox{\quad contains at least one prime.}
\end{equation}
\begin{equation}
\forall n \in[n_2,\infty), \quad \forall x > [90^{90}]^{n\over n-90}, \qquad 
[x,  x+ \; x^{1-{1\over n}}] \hbox{\quad contains at least one prime.}
\end{equation}
Note these are all lower bounds on the locations of prime-containing short intervals.
To improve on these bounds one would either need additional explicit results on primes between powers, (that is, reduce the exponent 90 even further), or to somehow reduce the constant $\exp(\exp(33))$ appearing in Dudek's result for primes between cubes.

%\clearpage 
If one restricts attention to integer values of $n$, (there is no pressing need to do so), then as alluded to in the abstract, one obtains slightly tidier looking (but slightly weaker) results
\begin{equation}
\forall n \in[4,90], \quad \forall x > \exp(3\exp(33)), \qquad 
[x,  x+ \; x^{1-{1\over n}}] \hbox{\quad contains at least one prime.}
\end{equation}
\begin{equation}
\forall n \in[91,\infty), \quad \forall x > [90^{90}]^{n\over n-90}, \qquad 
[x,  x+ \; x^{1-{1\over n}}] \hbox{\quad contains at least one prime.}
\end{equation}
We shall now further improve these bounds by using extra information coming from maximal prime gaps and effective linear intervals. 

%\clearpage
%================================================
\section{Use of maximal prime gaps}
%================================================

\enlargethispage{30pt}
If the interval $[x,  x+ \; x^{1-{1\over n}}] $ is asserted to contain at least one prime, then by setting $x=p_i+\epsilon$ for $0<\epsilon\ll1$ we see that the prime gaps satisfy $g_i \leq (p_i+\epsilon)^{1-{1\over n}}$, for all such~$\epsilon$. Consequently  $g_i \leq p_i^{1-{1\over n}}$. 

For current purposes we might as well set $n=2$ and consider $g_i \leq \sqrt{p_i}$ since larger exponents $n$ provide weaker constraints that are easier to satisfy.
But the (partial) truth of this statement can easily be verified by inspecting a table of maximal prime gaps to verify that it holds up to the location of the largest known maximal prime gap. 

The current (July 2025)  largest known maximal prime gap is the 85th maximal prime gap, located at $p_{85}^*= 101,412,319,996,363,309,069 \approx  10^{20}  \gtrsim 2^{66}$, with $n^*_{85}= \pi(p^*_{85})=2,251,483,061,895,611,799
\approx 2.25\times 10^{18}\gtrsim 2^{60}$, and gap size $g_{85}^*=1,854$. 
See references~\cite{prime-gaps, tables, primecount}. 

One wishes to check that 
\begin{equation}
\forall x < p_{85}^*, \qquad 
[x,  x+ \; x^{{1\over 2}}] \hbox{\quad contains at least one prime.}
\end{equation}
To do this, inspect the entire table of known maximal prime gaps, and simply verify that 
\begin{equation}
\forall i \in [1,85]  \qquad (p_i^*)^{{1\over 2}} > g_i^*.
\end{equation}
This is easily checked to be true.

Consequently we certainly have 
\begin{equation}
\forall n \geq 2, \quad \forall x < 10^{20}, \qquad 
[x,  x+ \; x^{1-{1\over n}}] \hbox{\quad contains at least one prime.}
\end{equation}
Note this is now an effective upper bound on the location of prime-containing short intervals.

%================================================
\section{Use of effective linear intervals}
%================================================
To try to bridge the gap between $10^{20}$ and the region above $90^{90}$ a viable strategy is to use high-precision linear intervals of the form $[(1-\Delta^{-1}) x,x]$.
Suitable intervals of this form have been studied by Cully-Hugill and Lee~\cite{Cully-Hugill+Lee}.

Assume we have data $(x_\Delta,\Delta)$ showing
\begin{equation}
[(1-\Delta^{-1}) x,x] \hbox{ contains a prime for } x > x_\Delta.
\end{equation}
Note this implies
\begin{equation}
\left[x,{x\over (1-\Delta^{-1})} \right] \hbox{ contains a prime for } x > x_\Delta (1-\Delta^{-1}).
\end{equation}
Thence
\begin{equation}
\left[x,x + {1\over \Delta-1} x  \right] \hbox{ contains a prime for } x > x_\Delta (1-\Delta^{-1}).
\end{equation}
So certainly
\begin{equation}
\left[x,x + {x\over\Delta-1} \right] \hbox{ contains a prime for } x > x_\Delta.
\end{equation}

\clearpage
For \emph{small} $x$ the Hoheisel interval $[x,x+x^{1-{1\over n}}]$ is ultimately 
larger  than the linear interval 
$\left[x,x + {x\over\Delta-1} \right] $. Equality occurs when
\begin{equation}
x^{1-{1\over n}} = {x\over\Delta-1}.
\end{equation}
That is, when
\begin{equation}
x^{1\over n} = \Delta-1
\end{equation}
implying
\begin{equation}
x= (\Delta-1)^{n}
\end{equation}
So given the $(x_\Delta,\Delta)$  data, the interval $[x,x+x^{1-{1\over n}}]$ certainly contains a prime in the (possibly empty) region
\begin{equation}
x \in [x_\Delta, (\Delta-1)^{n}] .
\end{equation}
A suitable collection of pairs $(x_\Delta,\Delta)$ is reported in reference~\cite{Cully-Hugill+Lee}. See Table 2. 

\begin{table}[!htbp]
\caption{Suitable set of pairs $(x_\Delta,\Delta)$. See reference~\cite{Cully-Hugill+Lee}.}
\begin{center}
\begin{tabular}{|c|c|}
\hline\hline
$\ln(x_\Delta)$ & $\Delta$\\
\hline\hline
$\ln(4\cdot 10^{18})$ &  $3.90970\cdot 10^7$\\
43& $  4.18168\cdot 10^7 $\\
46& $  1.63940\cdot 10^8 $\\
50& $  1.06120\cdot 10^9 $\\
55& $  1.02884\cdot 10^{10} $\\
60& $  7.69184\cdot 10^{10} $\\
75& $  1.74043\cdot 10^{11} $\\
90& $  1.84304\cdot 10^{11} $\\
105& $  1.91886\cdot 10^{11} $\\
120& $  1.97917\cdot 10^{11} $\\
135& $  2.02553\cdot 10^{11} $\\
150& $  2.07053\cdot 10^{11} $\\
300& $  2.30126\cdot 10^{11} $\\
600& $  2.51949\cdot 10^{11} $\\
\hline\hline
\end{tabular}
\end{center}
\label{default}
\end{table}%

To apply the $(x_\Delta,\Delta)$ data one picks a value of $n$,
 calculates all of the intervals $[x_\Delta, (\Delta-1)^{n}]$ from Table 2, checking that they are non-empty up to some $\Delta_{max}$, and that they overlap. The union of these intervals is $[4\cdot 10^{18}, (\Delta_{max}-1)^n]$.
 
 But then the interval
 \begin{equation}
 [x,x+x^{1-{1\over n}}] \quad\hbox{contains a prime for} \quad x \in [4\cdot 10^{18}, (\Delta_{max}-1)^n].
 \end{equation}
 A few key values are reported in Table 3.

\begin{table}[!htbp]
\caption{A few key instances of the intermediate bound. }
\begin{center}
\begin{tabular}{|c|c|c|}
\hline\hline
$n$ & $(\Delta_{max}-1)^n$ & Notes\\
\hline\hline
90   &  $1.312615432\cdot10^{1026}$ & C-H+J exponent\\
91   &  $3.307121456\cdot10^{1037}$ & Lowest relevant exponent\\
100 & $1.352894389\cdot10^{1140}$ & --- \\
105 & $1.373495024\cdot10^{1197}$ & Transitional exponent\\
106 & $3.460506978\cdot10^{1208}$ & Transitional exponent\\
125 & $1.459082570\cdot10^{1425}$ & --- \\
250 & $2.128921946\cdot10^{2850}$ & Heilbronn exponent\\
33000 & $5.411673381\cdot10^{250540}$& Hoheisel exponent\\
$\infty$ &  $\infty$ & --- \\
\hline\hline
\end{tabular}
\end{center}
\label{T:2}
\end{table}%

\newpage
%================================================
\section{Combining all of the above}
%================================================

For $n=106$, combining all the information from (i) the prime gaps, (ii) the $(x_\Delta,\Delta)$ long intervals, and (iii) the lower bound in Table 1, one sees that the entire real line ($x\geq 1$) is covered and so verifies that
\begin{equation}
\hbox{For } n=106, \quad \forall x \geq 1 \qquad  [x,x+x^{1-{1\over n}}] \quad\hbox{contains a prime}.
\end{equation}
The way the argument is set up this also automatically works for any $n>106$
\begin{equation}
\hbox{For } n\geq 106, \quad \forall x \geq 1 \qquad  [x,x+x^{1-{1\over n}}] \quad\hbox{contains a prime}.
\end{equation}
In particular, both the Hoheisel and Heilbronn results hold for all $x\geq 1$.
This now provides a fully effective version of both  the Hoheisel and Heilbronn results. 

For $n=105$, combining the information from the maximal prime gaps and the $(x_\Delta,\Delta)$ effective linear intervals covers the interval $x\in [1, \; 1.373495024\cdot10^{1197}]$, whereas the  lower bound in Table 1 covers the interval $x\in [1.488610050\cdot10^{1231},\; \infty)$. 
This leaves a ``smallish'' interval $x \in (1.373495024\cdot10^{1197}, \; 1.488610050\cdot10^{1231})$ wherein the arguments of this present article are inconclusive. 
Finding a few more  $(x_\Delta,\Delta)$ effective linear intervals would seem to be the most direct  plausible strategy for tightening this result. (Alternatively one could more indirectly try to decrease the exponent 90 appearing in the primes between powers result~\cite{CHJ2},  or the $\exp(\exp(33))$ appearing in Dudek's primes between cubes theorem~\cite{Dudek}).

%\clearpage
%================================================
\section{Discussion}
%================================================

\enlargethispage{50pt}
The good news is that the analysis above has yielded several explicit and effective bounds on when  sub-linear ``short''  intervals of the form $[x,x+x^{1-{1\over n}}]$ are guaranteed to contain a prime. 
The not so good news is that to obtain a widely applicable result $(x\geq 1)$ one needs a rather high exponent ($n\geq 106$), while to obtain a small exponent (say $n=4$) one needs a rather high value of $x$, $(x\geq \exp(3\exp(33)) \,)$.
Still, any progress along these lines should be of considerable interest.

%=================================================
\bigskip
\hrule\hrule\hrule
%================================================

\appendix
%================================================
\section{Primes between powers redux}
%================================================

Let us return to the issue of primes between powers. 

First, working ``downwards'' from infinity, note that Dudek's result for cubes
\begin{equation}
\forall x\geq  \exp(3\exp(33))   \qquad [x,  x+ 3 \; x^{1-{1\over3}}] \hbox{\quad contains at least one prime,}
\end{equation}
automatically implies
\begin{equation}
\forall n \geq 3 \quad \forall x\geq  \exp(3\exp(33))   \qquad [x,  x+ n \; x^{1-{1\over n}}] \hbox{\quad contains at least one prime,}
\end{equation}
whence
\begin{equation}
\forall n \geq 3 \quad \forall i\geq  \exp(3\exp(33)/n)   \qquad [i^n,  (i+1)^n] \hbox{\quad contains at least one prime.}
\end{equation}

We shall now seek to develop some additional partial results, working ``upwards'' from unity, by using the maximal prime gaps and the effective linear intervals of Table 2~\cite{Cully-Hugill+Lee}.
From the maximal prime gaps we easily see that
\begin{equation}
[i^2, (i+1)^2] \hbox{ contains a prime for } i \leq 10^{10} < \sqrt{p_{85}^*}.
\end{equation}
Thence certainly:
\begin{equation}
\hbox{For $n\geq 2$ we see that } [i^n, (i+1)^n] \hbox{ contains a prime for } i \leq  10^{20/n}.
\end{equation}
Regarding the effecitive linear intervals, we have already seen that if
\begin{equation}
[(1-\Delta^{-1}) x,x] \hbox{ contains a prime for } x > x_\Delta,
\end{equation}
then 
\begin{equation}
\left[x,x + {x\over\Delta-1} \right] \hbox{ contains a prime for } x > x_\Delta.
\end{equation}

For \emph{small} $x$ the primes-between-powers interval $[x,x+n x^{1-{1\over n}}]$ is ultimately 
larger  than the linear interval 
$\left[x,x + {x\over\Delta-1} \right] $. Equality occurs when
\begin{equation}
n x^{1-{1\over n}} = {x\over\Delta-1}.
\end{equation}
That is, when
\begin{equation}
x^{1\over n} = {n(\Delta-1)},
\end{equation}
implying
\begin{equation}
x= \left([\Delta-1] n\right)^{n}
\end{equation}
So given the $(x_\Delta,\Delta)$  data of Table 2~\cite{Cully-Hugill+Lee}, the interval $[x,x+n \,x^{1-{1\over n}}]$ certainly contains a prime in the (possibly empty) region
\begin{equation}
x \in \left[x_\Delta, \left([\Delta-1] n\right)^{n}\right] .
\end{equation}
This implies that  the 
 interval $[i^n,(i+1)^n]$ certainly contains a prime in the (possibly empty) region
\begin{equation}
i \in \left[\sqrt[n]{x_\Delta}, n\left(\Delta-1\right)\right] .
\end{equation}
For each specific value of $n$ we now work our way through Table 2. calculating the intervals  $\left[\sqrt[n]{x_\Delta}, n\left(\Delta-1\right)\right]$, checking that they are non empty, and that they overlap, out to some maximal value $\Delta_{max}$. 
Then the 
 interval $[i^n,(i+1)^n]$ certainly contains a prime for $i$ in the region
\begin{equation}
i \in \left[\sqrt[n]{10^{20}}, n\left(\Delta_{max}-1\right)\right] .
\end{equation}
For $n=2$ it turns out that this interval is empty and we get no extra information. 
For $n=3 $ this interval is not empty, so we do get some extra information, but a recent article~\cite{semiprimes-between-cubes} obtains somewhat better results.
For $n\geq 4 $ this interval is not empty and we do get extra information.

Combined with the result coming from the maximal prime gaps we see that for $n\geq 3$ the 
interval $[i^n,(i+1)^n]$ certainly contains a prime for $i$ in the region
\begin{equation}
i \in \left[1, n\left(\Delta_{max}-1\right)\right] .
\end{equation}
Specific results are reported in Table 4. 
Denoting $i_0(n) = n\left(\Delta_{max}-1\right)$ we note that $i_0(n)$ grows relatively slowly, while $i_0(n)^n$ grows explosively.

\bigskip
\hrule\hrule\hrule

\begin{table}[!htbp]
\caption{Regions, working ``upwards'' from unity, $i\in[1, i_0(n)]$ for which we are guaranteed the existence of primes between powers $[i^n,(i+1)^n]$.
For $n\geq 21$ we reach the end of the effective linear intervals reported in Table 2. 
For $n\geq 90$ the Cully-Hugill  and Johnston result~\cite{CHJ2} guarantees the universal existence of primes between powers $[i^n,(i+1)^n]$.}

\begin{center}
\begin{tabular}{|c|c|}
\hline\hline
$n$ & $ i_0(n) $ \\
\hline\hline
2   &  $10^{10}$ \\
\hline\hline
3 & $5.221290000\cdot10^{11}$ \\
\hline\hline
4 & $7.675440000\cdot10^{11}$ \\
5 & $1.012765000\cdot10^{12}$ \\
6 & $1.242318000\cdot10^{12}$  \\
7 & $1.449371000\cdot10^{12}$ \\
8& $1.656424000\cdot10^{12}$\\
9 &   $   1.863477000\cdot10^{12}$  \\
10 & $   2.070530000\cdot 10^{12} $\\
11 & $   2.531386000\cdot 10^{12} $\\
12& $    2.761512000\cdot 10^{12} $\\
13 & $   2.991638000\cdot 10^{12} $\\
14 & $   3.221764000\cdot 10^{12} $\\
15 & $   3.451890000\cdot 10^{12} $\\
16 & $   3.682016000\cdot 10^{12} $\\
17 & $   3.912142000 \cdot 10^{12} $\\
18 & $   4.142268000 \cdot 10^{12} $\\
19 & $   4.372394000 \cdot 10^{12} $\\
20 & $   4.602520000 \cdot 10^{12} $\\
%22 & $   5.542878000 \cdot 10^{12} $\\
\hline\hline
$n\in [21,89] $ & $ n \times  2.519490000 \cdot 10^{11} $\\
\hline\hline
21 & $   5.290929000 \cdot 10^{12} $\\
30 & $  7.558470000  \cdot 10^{12} $\\
40 & $  1.007796000  \cdot 10^{13} $\\
50 & $   1.259745000 \cdot 10^{13} $\\
60 & $   1.511694000 \cdot 10^{13} $\\
70 & $   1.763643000 \cdot 10^{13} $\\
80 & $    2.015592000\cdot 10^{13} $\\
89 & $    2.242346100\cdot 10^{13} $\\
\hline\hline
90 & $  \infty $\\
\hline\hline
\end{tabular}
\end{center}
\label{T:4}
\end{table}%

\vspace{5cm}

\clearpage

\null\vspace{-50pt}
%=================================================
\bigskip
\hrule\hrule\hrule
%================================================

%================================================
\section*{Acknowledgement}
%================================================
I wish to thank Daniel R. Johnston for bringing the updated results in references~\cite{CHJ1,CHJ2} to my attention.

\bigskip
%=================================================
\bigskip
\hrule\hrule\hrule
%================================================
%================================================
%\clearpage
%================================================

%================================================

\begin{thebibliography}{99}
%================================================
\newcommand{\arXiv}[1]{arXiv:~{\href{https://arxiv.org/abs/#1}{\color{blue}#1}}}
%This allows using \arXiv{2006.07125} or \arXiv{gr-qc/0009013} for nice links.
%================================================

\bibitem{Hoheisel}
Guido Hoheisel, 
``Primzahlprobleme in der Analysis'',\\
Sitzungsberichte Berliner Akad. d. Wiss., (1930) 580-588.

\bibitem{Heilbronn}
 Hans Heilbronn, ``Uber den Primzahlsatz von Herrn Hoheisel'',\\
 Math. Z. {\bf 36} (1933), 394--423.
\doi{10.1007/BF01188631}.

\bibitem{BHP}
R. C. Baker, G. Harman, and J. Pintz,\\
``The Difference Between Consecutive Primes, II",\\
Proc. London Math. Soc. {\bf 83 \# 3} (2001) 532--562.

\bibitem{Cully-Hugill}
Michaela Cully-Hugill,
``Primes between consecutive powers'',\\
Journal of Number Theory, {\bf 247} (2023) 100--117.
\doi{10.1016/j.jnt.2022.12.002}
[\arXiv{2107.14468} [math.NT]]

\bibitem{CHJ1}
Michaela Cully-Hugill and Daniel R. Johnston,\\
``On the error term in the explicit formula of Riemann--von Mangoldt'',\\
 International Journal of Number Theory
{\bf 19 \# 6} (2023) 1205--1228
\doi{10.1142/S1793042123500598}
[\arXiv{2111.10001} [math.NT]]

\bibitem{CHJ2}
Michaela Cully-Hugill and Daniel R. Johnston,\\
``On the error term in the explicit formula of Riemann--von Mangoldt. II.'', \\
Funct. Approx. Comment. Math. 
 {\bf 73 \# 2} (2025) 223--242 . 
\doi{10.7169/facm/241110-18-11}
[\arXiv{2402.04272} [math.NT]]


\bibitem{Dudek}
Adrian Dudek,
``An Explicit Result for Primes Between Cubes'',\\
 Funct. Approx. Comment. Math. {\bf 55 \# 2} (2016) 177--197. 
\doi{10.7169/facm/2016.55.2.3}
[\arXiv{1401.4233} [math.NT]]


\bibitem{prime-gaps}
Wikipedia, ``Prime gap'', 
\href{https://en.wikipedia.org/wiki/Prime_gap}{https://en.wikipedia.org/wiki/Prime\_gap} \\
(accessed on 12 July 2026). 

\bibitem{tables}
Prime pages, Table of known maximal gaps, 
\href{https://t5k.org/notes/GapsTable.html}{https://t5k.org/notes/GapsTable.html}\\
 (accessed on 12 July 2026). 
 
 \bibitem{primecount}
 Kim Walisch, 2025,  {\tt primecount} software,
 \href{https://github.com/kimwalisch/primecount}{https://github.com/kimwalisch/primecount}
 
 \bibitem{Cully-Hugill+Lee}
Michaela Cully-Hugill and Ethan S. Lee,\\
``Explicit Interval Estimates for Prime Numbers'',\\
Math. Comp. {\bf 91} (2022), 1955--1970.
\doi{10.1090/mcom/3719}\\{}
[\arXiv{2103.05986} [math.NT]]


\bibitem{semiprimes-between-cubes}
Daniel R. Johnston, Jonathan P. Sorenson, Simon N. Thomas, Jonathan E. Webster,
``Primes and almost primes between cubes'',
\arXiv{2601.15564} [math.NT]


%=================================================
\bigskip
\hrule\hrule\hrule
%================================================

%================================================
\end{thebibliography}
\end{document}